\documentclass{mfat}
\pagespan{295}{310}

\usepackage{mmap}
\usepackage[utf8]{inputenc}

\newtheorem{theorem}{Theorem}
\newtheorem{corollary}{Corollary}

\newtheorem{proposition}{Proposition}

\theoremstyle{remark}
\newtheorem{remark}{Remark}

\theoremstyle{definition}
\newtheorem{example}{Example}

\renewcommand{\rho}{\varrho}

\begin{document}

\title[Elliptic boundary-value problems in H\"ormander spaces]
{Elliptic boundary-value problems\\in H\"ormander spaces}


\author[A. Anop]{Anna Anop}
\address{Institute of Mathematics, National Academy of Sciences of Ukraine,
3 Tereshchenkivs'ka, Kyiv, 01601, Ukraine; Chernihiv National
Pedagogical University, 53 Het'mana Polubotka, Chernihiv, 14013,
Ukraine}
\email{ahlv@ukr.net}


\author[T. Kasirenko]{Tetiana Kasirenko}
\address{Institute of Mathematics, National Academy of Sciences of Ukraine,
3 Tereshchenkivs'ka, Kyiv, 01601, Ukraine}
\email{tania.nemtseva@gmail.com}


\date{19/09/2016}



\begin{abstract}
We investigate general elliptic boundary-value problems in H\"ormander inner product spaces that form the extended Sobolev scale. The latter  consists of all Hilbert spaces that are interpolation spaces with
respect to the Sobolev Hilbert scale. We prove that the operator corresponding to an arbitrary elliptic problem is Fredholm in appropriate couples of the H\"ormander spaces and induces a collection of isomorphisms on the extended Sobolev scale. We obtain a local a priory estimate for generalized solutions to this problem and prove a theorem on their local regularity in the H\"ormander spaces. We find new sufficient conditions under which generalized derivatives (of a given order) of the solutions are continuous.
\end{abstract}

\maketitle
\smallskip

\section{Introduction}\label{1sec1}

A fundamental result of the theory of general elliptic boundary-value problems consists in that these problems are Fredholm on appropriate couples of Sobolev spaces and induce isomorphisms between their subspaces of finite co-dimension (see monographs \cite{Berezansky68, Hermander63, Hermander85, LionsMagenes72, Triebel95, Roitberg96} and survey \cite{Agranovich97}). This result has important applications to the investigation of regularity of solutions to elliptic equations, to the spectral theory of elliptic differential operators and others, the most profound results being obtained for the Hilbert scale of Sobolev spaces. Nevertheless, the Sobolev scale is not finely calibrated enough for various problems appearing in the theory of partial differential equations and mathematical analysis (see monographs \cite{Hermander63, Hermander83, Jacob010205, MikhailetsMurach14, NicolaRodino10, Paneah00, Triebel01}).

In this connection, H\"ormander \cite{Hermander63, Hermander83} introduced and investigated a broad class of normed function spaces for which the index of regularity is a sufficiently arbitrary function depending on frequency variables. This function parameter allows one to characterize the regularity more finely than it is possible within the Sobolev space. H\"ormander gave applications of these spaces to the investigation of regularity of solutions to partial differential equations. But H\"ormander spaces have not found applications to boundary-value problems within decades. There was only one paper \cite{Slenzak74} by Slenzak devoted to the study of regular elliptic problems in some H\"ormander inner product spaces.

Recently Mikhailets and Murach \cite{MikhailetsMurach05UMJ5, MikhailetsMurach06UMJ3, MikhailetsMurach06UMJ11, MikhailetsMurach07UMJ5, MikhailetsMurach08UMJ4, MikhailetsMurach12BJMA2, MikhailetsMurach14} have built a theory of solvability of general elliptic boundary-value problems in a scale of H\"ormander inner product spaces whose regularity index is an arbitrary radial function varying regularly at infinity in the sense of Karamata. This refined Sobolev scale is obtained by the method of interpolation with a function parameter between Sobolev inner product spaces. This method plays a key role in the mentioned theory.

Therefore of definite interest is the class of all Hilbert spaces that are interpolation spaces between inner product Sobolev spaces. This class is constructively described by Mikhailets and Murach \cite{MikhailetsMurach09Dop3, MikhailetsMurach13UMJ3, MikhailetsMurach15ResMath1} with the help of Ovchinnikov's theorem \cite[Section~11.4]{Ovchinnikov84} and is called the extended Sobolev scale. This scale is formed by H\"ormander inner product spaces whose regularity index is an arbitrary radial function RO-varying at infinity in the sense of Avakumovi\'c. This scale contains the refined Sobolev scale and allows us to characterize the regularity of distributions more finely than the latter scale does. So, an RO-varying function need not have a number order of varying at infinity.

The purpose of the present paper is to investigate the character of solvability of general elliptic boundary-value problems and the properties of their generalized solutions on the extended Sobolev scale. Some results of the paper are announced in \cite{Anop14Dop4}. Various classes of elliptic problems were considered earlier in \cite{Anop13Coll2, AnopMurach14MFAT, AnopMurach14UMJ, Anop14Coll2} under stronger assumptions about the regularity index than the assumption made in the present paper. Note that elliptic systems given on $\mathbb{R}^{n}$ or a closed smooth manifold are investigated on this scale by Murach and Zinchenko \cite{MurachZinchenko13MFAT1, ZinchenkoMurach12UMJ11, ZinchenkoMurach14JMathSci}.

The paper consists of six sections. Section~\ref{1sec1} is Introduction. In Section~\ref{1sec2}, we state a general elliptic boundary-value problem and consider a formally adjoint problem with respect to the special Green formula. In section~\ref{1sec3}, we give a definition of H\"ormander spaces which form the extended Sobolev scale and discuss some of their properties. Section~\ref{1sec4} contains our main results concerning the properties of this problem on the extended Sobolev scale. These results are theorems on the Fredholm property of the problem and on isomorphisms induced by the problem, a local a priory estimate for its generalized solutions, and a theorem on their local regularity. As an application of H\"ormander spaces, we obtain new sufficient conditions under which generalized derivatives (of a prescribed order) of these solutions are continuous; specifically we give conditions for the generalized solutions to be classical. Section~\ref{1sec5} is devoted to the method of interpolation with a function parameter between Hilbert spaces. This method plays a key role in our investigation. The main results of the paper are proved in the last Section~\ref{1sec6}.

\section{Statement of the problem}\label{1sec2}

Let $\Omega\subset\mathbb{R}^n$ be a bounded domain with an infinitely smooth boundary $\Gamma$. (So, $\Gamma$ is a closed $C^{\infty}$-manifold of dimension $n-1$.) We suppose that the integer $n\geq2$. As usual, $\overline{\Omega}=\Omega\cup\Gamma$. Let $\nu(x)$ denote the unit vector of the inward normal to the boundary $\Gamma$ at a point $x\in\Gamma$.

We consider the following boundary-value problem:
\begin{gather}\label{1f1}
Au=f\quad\mbox{in}\quad\Omega,\\
B_{j}u=g_{j}\quad\mbox{on}\quad\Gamma,
\quad j=1,\ldots ,q.\label{1f2}
\end{gather}
Here,
$$
A:=A(x,D):=\sum_{|\mu|\leq 2q}a_{\mu}(x)D^{\mu}\
$$
is a linear differential operator on $\overline{\Omega}$ of an even order $2q\geq2$. Besides, every
$$
B_{j}:=B_{j}(x,D)=\sum_{|\mu|\leq m_{j}}b_{j,\mu}(x)D^{\mu}\
$$
is a linear boundary differential operator on $\Gamma$ of order $m_{j}\leq2q-1$. The coefficients of these operators are complex-valued infinitely smooth functions given on $\overline{\Omega}$ and $\Gamma$, respectively.

Here, we use the following standard designations:
$\mu:=(\mu_{1},\ldots,\mu_{n})$ is a multi-index with nonnegative integer components, $|\mu|:=\mu_{1}+\cdots+\mu_{n}$,
$D^{\mu}:=D_{1}^{\mu_{1}}\cdots D_{n}^{\mu_{n}}$, $D_{k}:=i\partial/\partial x_{k}$, $k=1,\ldots ,n$, where $i$ is imaginary unit and $x=(x_1,\ldots,x_n)$ is an arbitrary point in $\mathbb{R}^{n}$. We also put $D_{\nu}:=i\partial/\partial\nu(x)$.

In the paper, we suppose that the boundary-value problem \eqref{1f1}, \eqref{1f2} is elliptic in the domain $\Omega$. This means that the differential operator $A$ is properly elliptic on $\overline{\Omega}$ and that the system of boundary differential operators $B:=(B_{1},\ldots,B_{q})$ satisfies Lopatynsky condition with respect to $A$ on $\Gamma$ (see, e.g., survey \cite[Section~1.2]{Agranovich97}).

With the problem \eqref{1f1}, \eqref{1f2}, we associate the linear mapping
\begin{equation}\label{1f3}
\begin{gathered}
u\mapsto(Au,Bu)=(Au,B_{1}u,\ldots ,B_{q}u),\quad
\mbox{with}\quad u\in C^{\infty}(\overline{\Omega}).
\end{gathered}
\end{equation}
We investigate properties of the extension (by continuity) of this mapping in appropriate couples of H\"ormander inner product spaces which form the extended Sobolev scale.

In order to describe the range of this extension, we need the following special Green formula \cite[Theorem~3.1.1]{KozlovMazyaRossmann97}:
$$
(Au,v)_{\Omega}+\sum_{j=1}^{q}(B_{j}u,h_{j})_{\Gamma}=
(u,A^{+}v)_{\Omega}+\sum_{j=1}^{2q}\biggl(D_{\nu}^{j-1}u,K_{j}v+
\sum_{k=1}^{q}Q_{k,j}^{+}h_{k}\biggr)_{\Gamma},
$$
where $u,v\in C^{\infty}(\overline{\Omega})$, $h_{1},\ldots,h_{q}\in C^{\infty}(\Gamma)$, and, besides, $(\cdot,\cdot)_{\Omega}$ and $(\cdot,\cdot)_{\Gamma}$ denote the inner products in the Hilbert spaces $L_{2}(\Omega)$ and $L_{2}(\Gamma)$ of square integrable functions over $\Omega$ and $\Gamma$ respectively and later denote the extension of these inner products by continuity. As usual, the differential operator $A^{+}$ is formally adjoint to $A$; namely,
$$
(A^{+}v)(x):=
\sum_{|\mu|\leq2q}D^{\mu}\bigl(\overline{a_{\mu}(x)}v(x)\bigr).
$$
Besides, each tangent differential operator $Q_{k,j}^{+}$ is adjoint to the tangent differential operator $Q_{k,j}:=Q_{j,k}(x,D_{\tau})$ appearing in the representation
$$
B_{j}(x,D)=
\sum_{k=1}^{2q}Q_{j,k}(x,D_{\tau})D_{\nu}^{k-1}.
$$
(Of course, if $k>m_{j}+1$, then $Q_{j,k}:=0$.) Furthermore, each $K_{j}:=K_{j}(x,D)$ is a certain linear boundary differential operator of order $2q-j$ with coefficients from $C^{\infty}(\Gamma)$.

Note that we do not suppose the system $B$ of boundary operators to be normal. Therefore the classical Green formula \cite[Section~4.2]{Agranovich97} is not generally valid for the problem \eqref{1f1}, \eqref{1f2}.

In view of the special Green formula, we consider the following boundary-value problem in $\Omega$ with $q$ additional unknown functions on $\Gamma$:
\begin{gather}\label{1f4}
A^{+}v=\omega\quad\mbox{in}\quad\Omega,\\
\label{1f5}
K_{j}v+\sum_{k=1}^{q}Q_{k,j}^{+}h_{k}=\zeta_{j}\quad
\mbox{on}\quad\Gamma,\quad j=1,\ldots ,2q.
\end{gather}
This problem is formally adjoint to the problem \eqref{1f1}, \eqref{1f2} with respect to the given Green formula. The boundary-value problem \eqref{1f1}, \eqref{1f2} is elliptic if and only if the formally adjoint problem \eqref{1f4}, \eqref{1f5} is elliptic in the relevant sense \cite[Theorem~3.1.2]{KozlovMazyaRossmann97}.

\section{The extended Sobolev scale}\label{1sec3}

As we have noted, this scale was selected and investigated by Mikhailets and Murach \cite{MikhailetsMurach09Dop3, MikhailetsMurach13UMJ3, MikhailetsMurach14, MikhailetsMurach15ResMath1}. The scale consists of isotropic H\"ormander inner product spaces $H^{\alpha}$ for which the index $\alpha$ of regularity runs through the function class $\mathrm{RO}$.

By definition, the class $\mathrm{RO}$ consists of all Borel measurable functions $\alpha:\nobreak[1,\infty)\rightarrow(0,\infty)$ for which there exist numbers $b>1$ and $c\geq1$ such that
$$
c^{-1}\leq\frac{\alpha(\lambda t)}{\alpha(t)}\leq c\quad\mbox{for every}\quad t\geq1\quad\mbox{and}\quad\lambda\in[1,b]
$$
($b$ and $c$ may depend on $\alpha$). Such functions are called RO-varying (or OR-varying) at infinity. The class $\mathrm{RO}$ was introduced by
Avakumovi\'c \cite{Avakumovic36} in 1936 and is sufficiently
investigated (see, e.g., monographs \cite[Sections 2.0--2.2]{BinghamGoldieTeugels89} and \cite[Appendix~1]{Seneta76}).

This class admits a simple description \cite[Theorem~A.1]{Seneta76}, namely
$$
\alpha\in\mathrm{RO}\;\;\Leftrightarrow\;\;\alpha(t)=\exp\Biggl(\beta(t)+
\int_{1}^{\:t}\frac{\gamma(\tau)}{\tau}\;d\tau\Biggr)\quad
\mbox{for}\quad t\geq1,
$$
where the real-valued functions $\beta$ and $\gamma$ are Borel measurable and bounded on $[1,\infty)$.

We need the following property of the class $\mathrm{RO}$ \cite[Theorem~A.2]{Seneta76}: for every function $\alpha\in\mathrm{RO}$ there exist numbers $s_{0},s_{1}\in\mathbb{R}$,
$s_{0}\leq s_{1}$, and $c_{0},c_{1}>0$ such that
\begin{equation}\label{Ax=b3}
c_{0}\lambda^{s_{0}}\leq\frac{\alpha(\lambda t)}{\alpha (t)}\leq
c_{1}\lambda^{s_{1}} \quad\mbox{for all}\quad t\geq1,\quad \lambda\geq 1.
\end{equation}
We put
\begin{gather*}
\sigma_{0}(\alpha):=
\sup\,\{s_{0}\in\mathbb{R}:\,\mbox{the left-hand inequality in \eqref{Ax=b3} holds}\},\\
\sigma_{1}(\alpha):=\inf\,\{s_{1}\in\mathbb{R}:\,\mbox{the right-hand inequality in
\eqref{Ax=b3} holds}\}.
\end{gather*}
The numbers
$\sigma_{0}(\alpha)$ and $\sigma_{1}(\alpha)$ are the lower and upper Matuszewska indices of the function $\alpha\in\mathrm{RO}$ (see \cite{Matuszewska64} and \cite[Section~2.1.2]{BinghamGoldieTeugels89}). Of course, $-\infty<\sigma_{0}(\alpha)\leq\sigma_{1}(\alpha)<\infty$.

We consider some examples of RO-varying functions.

\begin{example}\label{ex1} \rm
Let a continuous function $\alpha:[1,\infty)\rightarrow(0,\infty)$ be such that
\begin{equation*}
\alpha(t):=t^{s}(\log t)^{r_{1}}(\log\log
t)^{r_{2}}\ldots(\underbrace{\log\ldots\log}_{k\;\mathrm{times}} t)^{r_{k}}\quad\mbox{for}\quad t\gg1.
\end{equation*}
Here, the integer $k\geq1$ and real numbers $s,r_{1},\ldots,r_{k}$ are arbitrarily chosen. Then $\alpha\in\mathrm{RO}$ and $\sigma_{0}(\alpha)=\sigma_{1}(\alpha)=s$.
\end{example}

\begin{example}\label{ex2} \rm
Let $\theta\in\mathbb{R}$, $\delta>0$ and $r\in(0,1]$. We put
\begin{equation*}
\alpha(t):=\left\{
\begin{array}{ll}
t^{\theta+\delta\sin(\log\log t)^{r}}\; &\hbox{for}\;\;t>e,\\
t^{\theta}\; &\hbox{for}\;\;1\leq t\leq e.
\end{array}\right.
\end{equation*}
Then $\alpha\in\mathrm{RO}$, with $\sigma_{0}(\alpha)=\theta-\delta$ and $\sigma_{1}(\alpha)=\theta+\delta$.
\end{example}

Let $\alpha\in\mathrm{RO}$. We will first introduce the H\"ormander space $H^{\alpha}$ over $\mathbb{R}^{n}$, where $1\leq n\in\mathbb{Z}$, and then over $\Omega$ and $\Gamma$.

By definition, the complex linear space
$H^{\alpha}(\mathbb{R}^{n})$ consists of all tempered distributions in $\mathbb{R}^{n}$ such that their Fourier transform $\widehat{w}$
is locally Lebesgue integrable over $\mathbb{R}^{n}$ and satisfies the condition
$$
\|w\|_{H^{\alpha}(\mathbb{R}^{n})}^{2}:=
\int_{\mathbb{R}^{n}}\alpha^2((1+|\xi|^{2})^{1/2})\,
|\widehat{w}(\xi)|^2\,d\xi
<\infty.
$$
This space is Hilbert and separable with respect to the norm $\|\cdot\|_{H^{\alpha}(\mathbb{R}^{n})}$.
In the paper we interpret  distributions as \emph{anti}linear functionals given on the relevant  spaces of test functions.

The space $H^{\alpha}(\mathbb{R}^{n})$ is a Hilbert and isotropic case of the spaces $B_{p,k}$ introduced and
investigated by H\"ormander
(see his monographs \cite[Sec.~2.2]{Hermander63} and \cite[Sec.~10.1]{Hermander83}). Namely, $H^{\alpha}(\mathbb{R}^{n})=B_{p,k}$ if $p=2$ and $k(\xi)=\alpha((1+|\xi|^{2})^{1/2})$ for all $\xi\in\mathbb{R}^{n}$. Remark that, in the Hilbert case of $p=2$, the H\"ormander spaces coincide with the spaces introduced and investigated by Volevich and Paneah \cite[Sec.~2]{VolevichPaneah65}.

If $\alpha(t)\equiv t^{s}$ for some $s\in\mathbb{R}$, then  $H^{\alpha}(\mathbb{R}^{n})$ becomes the inner product Sobolev space $H^{(s)}(\mathbb{R}^{n})$ of order $s$. Generally, \vspace*{-1mm}
\begin{equation}\label{1f7}
(s_{0}<\sigma_{0}(\alpha),\;\;\sigma_{1}(\alpha)<s_{1})\;\Rightarrow\; H^{(s_{1})}(\mathbb{R}^{n})\hookrightarrow H^{\alpha}(\mathbb{R}^{n})\hookrightarrow H^{(s_{0})}(\mathbb{R}^{n}),
\end{equation}
with embeddings being continuous and dense.

Consider the class of Hilbert function spaces $\{H^{\alpha}(\mathbb{R}^{n}):\alpha\in\mathrm{RO}\bigr\}$.
Following Mik\-hai\-lets and Murach \cite{MikhailetsMurach13UMJ3, MikhailetsMurach14}, we call it the extended Sobolev scale over $\mathbb{R}^{n}$. Its analogs for the Euclidean domain $\Omega$ and the closed infinitely smooth manifold $\Gamma$ are introduced in the standard way; see \cite[Sec.~2]{MikhailetsMurach15ResMath1} and \cite[Sec.~2.4.2]{MikhailetsMurach14} respectively). Let us give the necessary definitions.

By definition, the linear space $H^{\alpha}(\Omega)$ consists of the restrictions of all distributions $w\in H^{\alpha}(\mathbb{R}^{n})$ to the domain $\Omega$. The norm in $H^{\alpha}(\Omega)$ is defined by the  formula
$$
\|v\|_{H^{\alpha}(\Omega)}:=
\inf\bigl\{\,\|w\|_{H^{\alpha}(\mathbb{R}^{n})}:\,
w\in H^{\alpha}(\mathbb{R}^{n}),\ w=v\;\,\mbox{in}\;\,\Omega\,\bigr\},
$$
with $v\in H^{\alpha}(\Omega)$. The space $H^{\alpha}(\Omega)$
is Hilbert and separable with respect to this norm; the set $C^{\infty}(\overline{\Omega})$ is dense in $H^{\alpha}(\Omega)$.

In short, the linear space $H^{\alpha}(\Gamma)$ consists of all distributions on $\Gamma$ that belong to $H^{\alpha}(\mathbb{R}^{n-1})$ in local coordinates. Let us give the detailed definition. We arbitrarily choose a finite atlas from $C^{\infty}$-structure on  $\Gamma$ formed by local charts $\pi_j: \mathbb{R}^{n-1}\leftrightarrow \Gamma_{j}$, where
$j=1,\ldots,\varkappa$. Here the open sets $\{\Gamma_{1},\ldots,\Gamma_{\varkappa}\}$ form a covering of the manifold $\Gamma$. Let functions $\chi_j\in C^{\infty}(\Gamma)$, with $j=1,\ldots,\varkappa$, form a partition of unity on $\Gamma$ that satisfies the condition $\mathrm{supp}\,\chi_j\subset \Gamma_j$ for each $j\in\{1,\ldots,\varkappa\}$.

By definition, the complex linear space $H^{\alpha}(\Gamma)$ consists of all distributions $h$ on $\Gamma$ such that $(\chi_{j}h)\circ\pi_{j}\in H^{\alpha}(\mathbb{R}^{n-1})$ for all
$j\in\{1,\ldots,\varkappa\}$. Here, $(\chi_{j}h)\circ\pi_{j}$ is the representation of the distribution $\chi_{j}h$ in the local chart $\pi_{j}$. The space $H^{\alpha}(\Gamma)$ is endowed with the norm
$$
\|h\|_{H^{\alpha}(\Gamma)}:=
\biggl(\,\sum_{j=1}^{\varkappa}\,
\|(\chi_{j}h)\circ\pi_{j}\|_
{H^{\alpha}(\mathbb{R}^{n-1})}^{2}\biggr)^{1/2}.
$$
This space is Hilbert and separable with respect to this norm and does not depend (up to equivalence of norms) on our choice of the atlas and the partition of unity \cite[Theorem~2.21]{MikhailetsMurach14}. The set $C^{\infty}(\Gamma)$ is dense in $H^{\alpha}(\Gamma)$.

The above-defined function spaces form the extended Sobolev scales
$$
\{H^{\alpha}(\Omega):\alpha\in\mathrm{RO}\}\quad\mbox{and}\quad
\{H^{\alpha}(\Gamma):\alpha\in\mathrm{RO}\}
$$
over $\Omega$ and $\Gamma$ respectively. They contain the Hilbert scale of Sobolev space: if $\alpha(t)\equiv t^{s}$ then $H^{\alpha}(\Omega)=:H^{(s)}(\Omega)$ and
$H^{\alpha}(\Gamma)=:H^{(s)}(\Gamma)$ are the inner product Sobolev spaces over $\Omega$ or $\Gamma$ of order $s\in\mathbb{R}$.

Dealing with the extended Sobolev scale, we will use some of its properties relating to embeddings of spaces. Let $\alpha,\eta\in\mathrm{RO}$ and $G\in\{\Omega,\Gamma\}$. The function $\alpha/\eta$ is bounded in a neighborhood of infinity if and only if $H^{\eta}(G)\hookrightarrow H^\alpha(G)$. This embedding is continuous  and dense. It is compact if and only if
$\alpha(t)/\eta(t)\rightarrow0$ as $t\rightarrow\infty$.
This property is a consequence of Theorems 2.2.2 and 2.2.3 from H\"ormander's monograph \cite{Hermander63}. Specifically, property \eqref{1f7} remains true if we replace $\mathbb{R}^{n}$ with $\Omega$ or $\Gamma$. Then the embeddings are compact and dense.

We also mention H\"ormander's result \cite[Theorem~2.2.7]{Hermander63} on the embeddings of his spaces in spaces of continuously differentiable functions. For the extended Sobolev scale over $\Omega$, this result is formulated as follows:
\begin{equation}\label{v}
\int_1^{\infty} t^{2p+n-1}\alpha^{-2}(t)\,dt<\infty\;\;\Leftrightarrow\;\;
H^\alpha(\Omega)\hookrightarrow C^p(\overline{\Omega});
\end{equation}
here, $0\leq p\in\mathbb{Z}$ and $\alpha\in\mathrm{RO}$. This embedding is continuous (see \cite[Proposition~2.6(vi)]{MikhailetsMurach14} and \cite[Lemma~2]{ZinchenkoMurach12UMJ11}). Note that, in the Sobolev case where $\alpha(t)\equiv t^{s}$ for some $s\in\mathbb{R}$, property \eqref{v} becomes the well-know Sobolev embedding theorem according to which the inequality $s>p+n/2$ is equivalent to the embedding $H^{(s)}(\Omega)\hookrightarrow C^p(\overline{\Omega})$.

\section{Main results}\label{1sec4}

In this section we will formulate our main results concerning the properties of the elliptic boundary-value problem \eqref{1f1}, \eqref{1f2} on the extended Sobolev scale.

Let us previously agree about the following. We will use the H\"ormander spaces $H^{\alpha}$ for which the order of regularity is given in the form $\alpha(t)\equiv\varphi(t)t^{s}$ for some $\varphi\in\mathrm{RO}$ and $s\in\mathbb{R}$. To not write the argument $t$ in superscripts, we resort to the function parameter $\rho(t):=t$ of $t\geq1$. Then $H^{\alpha}$ can be written as $H^{\varphi\rho^{s}}$. Note that if $\varphi\in\mathrm{RO}$ and $s\in\mathbb{R}$, then $\varphi\rho^s\in\mathrm{RO}$ and
$\sigma_j(\rho^s\varphi)=s+\sigma_j(\varphi)$ for each $j\in\{0,1\}$.

We associate two finite-dimensional spaces with the boundary-value problem \eqref{1f1}, \eqref{1f2} and its adjoint problem \eqref{1f4}, \eqref{1f5}. Namely, let $N$ denote the linear space of all solutions
$u\in C^{\infty}(\overline{\Omega})$ to the problem \eqref{1f1}, \eqref{1f2} in the case where $f=0$ on $\Omega$ and each $g_{j}=0$ on~$\Gamma$. Besides, we let $N_{*}$ denote the linear space of all solutions
$$
(v,h_{1},\ldots,h_{q})\in C^{\infty}(\overline{\Omega})\times(C^{\infty}(\Gamma))^{q}
$$
to the adjoint problem \eqref{1f4}, \eqref{1f5} in the case where $\omega=0$ on $\Omega$ and each $\zeta_{j}=0$ on~$\Gamma$.
We put $m:=\mathrm{max}\{m_{1},\ldots ,m_{q}\}$; recall that each $m_{j}:=\mathrm{ord}\,B_{j}$. Since these problems are elliptic in $\Omega$, the spaces $N$ and $N_{*}$ are finite-dimensional (see, e.g., \cite[Lemma~3.4.2]{KozlovMazyaRossmann97}).

\begin{theorem}\label{1th1}
Let $\varphi\in\mathrm{RO}$ and
$\sigma_0(\varphi)>m+1/2$. Then the mapping \eqref{1f3} extends uniquely (by continuity) to a bounded operator
\begin{equation}\label{1f8}
(A,B):H^{\varphi}(\Omega)\rightarrow
H^{\varphi\rho^{-2q}}(\Omega)\oplus\bigoplus_{j=1}^{q}
H^{\varphi\rho^{-m_j-1/2}}(\Gamma)=:
\mathcal{H}^{\varphi\rho^{-2q}}(\Omega,\Gamma).
\end{equation}
This operator is Fredholm. Its kernel is equal to $N$, and its range consists of all vectors $(f,g_1,\ldots,g_q)\in\mathcal{H}^{\varphi\rho^{-2q}}(\Omega,\Gamma)$ such that
\begin{equation}\label{1f9}
(f,v)_\Omega+\sum_{j=1}^{q}(g_j,h_{j})_{\Gamma}=0 \quad
\mbox{for all} \quad (v,h_{1},\ldots,h_{q})\in N_{*}.
\end{equation}
The index of the operator \eqref{1f8} is equal to $\dim N-\dim N_{*}$ and does not depend on $\varphi$.
\end{theorem}

In view of this theorem, it is worthwhile to recall that a linear bounded operator $T:E_{1}\rightarrow E_{2}$, where $E_{1}$ and $E_{2}$ are Banach spaces, is called Fredholm if its kernel $\ker T$ and co-kernel $E_{2}/T(E_{1})$ are both finite-dimensional. If this operator is Fredholm, then its range is closed in $E_{2}$ (see, e.g., \cite[Lemma~19.1.1]{Hermander85}), and its index $\mathrm{ind}\,T:=\dim\ker T-\dim(E_{2}/T(E_{1}))$ is well defined and finite.

Note that the assumption $\sigma_0(\varphi)>m+1/2$ made in Theorem~\ref{1th1} cannot be weakened. Specifically, if $\varphi(t)\equiv t^{s}$ for some real $s\leq m+1/2$, then the mapping $u\mapsto Bu$, with $u\in C^{\infty}(\overline{\Omega})$, cannot be extended to a continuous linear operator from $H^{(s)}(\Omega)$ to $(\mathcal{D}'(\Gamma))^{q}$. Here, as usual, $\mathcal{D}'(\Gamma)$ denotes the linear topological space of all distribution on $\Gamma$.
Anop \cite{Anop13Coll2, AnopMurach14MFAT, AnopMurach14UMJ, Anop14Coll2} proved Theorem~\ref{1th1} for various classes of elliptic boundary-value problems under stronger restrictions on $\sigma_0(\varphi)$ than $\sigma_0(\varphi)>m+1/2$. At least, it was supposed that $\sigma_0(\varphi)>2q-1/2$.

Perhaps, the expression $(f,v)_\Omega$ in \eqref{1f9} needs commenting. If $\sigma_{0}(\varphi)>2q$, then $f\in H^{\varphi\varrho^{-2q}}(\Omega)\subset L_{2}(\Omega)$ and, hence, $(f,v)_\Omega$ is the inner product in $L_{2}(\Omega)$. If $m+1/2<\sigma_{0}(\varphi)\leq 2q$, then $(f,v)_\Omega$ is well defined by the passage to the limit; namely,  $(f,v)_\Omega:=\lim_{k\to\infty}(f_{k},v)_\Omega$ where $C^{\infty}(\overline{\Omega})\ni f_{k}\to f$ in $H^{\varphi\varrho^{-2q}}(\Omega)$ as $k\to\infty$. Here, the condition $(v,h_{1},\ldots,h_{q})\in N_{*}$ is essential. So, if $\sigma_{0}(\varphi)\leq 2q-1/2$ and $v$ is an arbitrarily chosen function from $C^{\infty}(\overline{\Omega})$, then the mapping  $f\mapsto(f,v)_\Omega$, with $f\in C^{\infty}(\overline{\Omega})$, cannot be extended to a continuous linear functional on $H^{\varphi\varrho^{-2q}}(\Omega)$. Specifically, for Sobolev spaces, this follows from \cite[Theorems 4.8.2(c) and 4.3.2/1(c)]{Triebel95}.

In the case where $N=\{0\}$ and $N_{*}=\{0\}$, the operator \eqref{1f8} is an isomorphism between spaces $H^{\varphi}(\Omega)$ and $\mathcal{H}^{\varphi\rho^{-2q}}(\Omega,\Gamma)$. This follows from  Theorem~\ref{1th1} and the Banach theorem on inverse operator. In the general case, the operator \eqref{1f8} induced an isomorphism between their certain subspaces of finite co-dimension. These subspaces and projectors onto them can be build in the following way.

Let $\varphi\in\mathrm{RO}$ and $\sigma_0(\varphi)>m+1/2$. We consider the decompositions of the spaces in the direct sums of their subspaces:
\begin{gather}\label{1f10}
H^{\varphi}(\Omega)=N\dotplus\bigl\{u\in
H^{\varphi}(\Omega):\,(u,v)_\Omega=0\;\;\mbox{for all}\;\;v\in N\bigr\},\\
\mathcal{H}^{\varphi\rho^{-2q}}(\Omega,\Gamma)=
N_{*}\dotplus\bigl\{(f,g_1,\ldots,g_q)\in
\mathcal{H}^{\varphi\rho^{-2q}}(\Omega,\Gamma):
\eqref{1f9}\;\mbox{is true}\bigr\}.\label{1f11}
\end{gather}
The representation \eqref{1f10} exists since it is a restriction of the decomposition of the space  $L_{2}(\Omega)$ into the orthogonal sum of the subspace $N$ and its complement. The equality \eqref{1f11} holds because the intersection of the subspaces on its right-hand side is zero space and because the finite dimension of the first subspace is equal to the co-dimension of the second due to Theorem~\ref{1th1}. Let $P$ and $P_{*}$ respectively denote the oblique projectors of the spaces
$H^{\varphi}(\Omega)$ and $\mathcal{H}^{\varphi\rho^{-2q}}(\Omega,\Gamma)$ onto the second summands in \eqref{1f10} and \eqref{1f11} parallel to the first. Evidently, these projectors do not depend on~$\varphi$.

\begin{theorem}\label{1th2}
For arbitrary $\varphi\in\mathrm{RO}$ with $\sigma_{0}(\varphi)>m+1/2$, the restriction of the mapping~\eqref{1f8} to the subspace
$P(H^{\varphi}(\Omega))$ is an isomorphism
\begin{equation}\label{1f12}
(A,B):\,P(H^{\varphi}(\Omega))\leftrightarrow
P_{*}\bigl(\mathcal{H}^{\varphi\rho^{-2q}}(\Omega,\Gamma)\bigr).
\end{equation}
\end{theorem}

Let us now consider properties of the generalized solutions to the elliptic boundary-value problem \eqref{1f1}, \eqref{1f2} on the extended Sobolev scale. We recall the definition of these solutions. We put
$$
H^{m+1/2+}(\Omega):=
\bigcup_{\substack{\alpha\in\mathrm{RO}:\\\sigma_{0}(\alpha)>m+1/2}}
H^{\alpha}(\Omega)=\bigcup_{s>m+1/2}H^{(s)}(\Omega)
$$
The latter equality is valid in view of \eqref{1f7}. Owing to Theorem~\ref{1th1}, for every function $u\in H^{m+1/2+}(\Omega)$ we can reasonably define the vector
\begin{equation*}
(f,g):=(f,g_{1},\ldots,g_{q}):=(A,B)u\in\mathcal{D}'(\Omega)\times
(\mathcal{D}'(\Gamma))^{q}.
\end{equation*}
Here, as usual, $\mathcal{D}'(\Omega)$ and $\mathcal{D}'(\Gamma)$ stand respectively for the linear topological spaces of all distributions on the domain $\Omega$ or on the manifold $\Gamma$. The function $u$ is called a generalized (strong) solution of the boundary-value problem \eqref{1f1}, \eqref{1f2} with the right-hand side $(f,g)$.

\begin{theorem}\label{1th3}
Let $\varphi\in\mathrm{RO}$ with $\sigma_{0}(\varphi)>m+1/2$, and let functions $\chi,\eta\in C^{\infty}(\overline{\Omega})$ satisfy the condition $\eta=1$ in a neighborhood of $\mathrm{supp}\,\chi$. Then there exists a number $c=c(\varphi,\chi,\eta)>0$ such that an arbitrary function
$u\in H^{\varphi}(\Omega)$ satisfies the estimate
\begin{equation}\label{1f14}
\|\chi u\|_{H^{\varphi}(\Omega)}\leq c\,\bigl(\|\chi(A, B)u\|_{\mathcal{H}^{\varphi\rho^{-2q}}(\Omega,
\Gamma)}+\|\eta u\|_{H^{\varphi\rho^{-1}}(\Omega)}\bigr).
\end{equation}
Here, $c$ does not depend on $u$.
\end{theorem}

In the special case where $\chi=\eta=1$, inequality \eqref{1f14} is a global a priori estimate of the generalized solution $u$ to the elliptic boundary-value problem \eqref{1f1}, \eqref{1f2}. Generally, this inequality is a local a priory estimate of the solution. Indeed, for every nonempty open subset of $\overline{\Omega}$, we can choose the functions $\chi,\eta$ so that they satisfy the condition of Theorem~\ref{1th3} and that their supports lie in this subset.

Let us discuss the local regularity of the generalized solutions. Let $V$ be an open subset of $\mathbb{R}^{n}$ that has a nonempty intersection with $\Omega$. We put $\Omega_0:=\Omega\cap V$ and $\Gamma_{0}:=\Gamma\cap V$ (the $\Gamma_{0}=\varnothing$ case is possible). For arbitrary parameter $\alpha\in\mathrm{RO}$, we introduce a local analogs of the spaces $H^{\alpha}(\Omega)$ and $H^{\alpha}(\Gamma)$ in the following standard way. By definition, the linear space $H^{\alpha}_{\mathrm{loc}}(\Omega_{0},\Gamma_{0})$ consists of all distributions $u\in\mathcal{D}'(\Omega)$ such that $\chi u\in H^{\alpha}(\Omega)$ for arbitrary function $\chi\in C^{\infty}(\overline{\Omega})$ with $\mathrm{supp}\,\chi\subset\Omega_0\cup\Gamma_{0}$. The topology in the linear space $H^{\alpha}_{\mathrm{loc}}(\Omega_{0},\Gamma_{0})$ is induced by the seminorms $u\mapsto\|\chi u\|_{H^{\alpha}(\Omega)}$, where $\chi$ is an arbitrary function from the definition of this space. Similarly, the linear space $H^{\alpha}_{\mathrm{loc}}(\Gamma_{0})$ consists of all distributions  $h\in\nobreak\mathcal{D}'(\Gamma)$ such that $\chi h\in H^{\alpha}(\Gamma)$ for arbitrary function $\chi\in C^{\infty}(\Gamma)$ with $\mathrm{supp}\,\chi\subset\Gamma_{0}$.
The topology in the linear space
$H^{\alpha}_{\mathrm{loc}}(\Gamma_{0})$ is induced by the seminorms $h\mapsto\|\chi h\|_{H^{\alpha}(\Omega)}$, where $\chi$ is an arbitrary function from the definition of this space.

\begin{theorem}\label{1th4}
Let a function $u\in H^{m+1/2+}(\Omega)$ be a generalized solution to the elliptic boundary-value problem \eqref{1f1}, \eqref{1f2} in which
$$
(f,g_1,\ldots,g_q)\in H^{\varphi\rho^{-2q}}_{\mathrm{loc}}(\Omega_{0},\Gamma_{0})
\oplus\bigoplus_{j=1}^{q}
H^{\varphi\rho^{-m_{j}-1/2}}_{\mathrm{loc}}(\Gamma_{0})=:
\mathcal{H}^{\varphi\rho^{-2q}}_{\mathrm{loc}}(\Omega_{0},\Gamma_{0})
$$
for a certain parameter $\varphi\in\mathrm{RO}$ such that $\sigma_{0}(\varphi)>m+1/2$. Then $u\in H^{\varphi}_{\mathrm{loc}}(\Omega_{0},\Gamma_{0})$.
\end{theorem}

It is worthwhile to mention two special cases of Theorem~\ref{1th4}. In the case where $\Omega_{0}=\Omega$ and $\Gamma_{0}=\Gamma$, the spaces $H^{\varphi}_{\mathrm{loc}}(\Omega_{0},\Gamma_{0})$ and $\mathcal{H}^{\varphi\rho^{-2q}}_{\mathrm{loc}}(\Omega_{0},\Gamma_{0})$ coincide respectively with $H^{\varphi}(\Omega)$ and $\mathcal{H}^{\varphi\rho^{-2q}}(\Omega,\Gamma)$. Therefore Theorem~\ref{1th4} asserts in this case that the regularity of $u$ increases globally, i.e. in the whole domain $\Omega$ up to its boundary $\Gamma$. The case where $\Omega_{0}=\Omega$ and $\Gamma_{0}=\varnothing$ leads us to the assertion that the smoothness of $u$ increases in neighborhoods of all internal points of $\overline{\Omega}$.

As an application of the extended Sobolev scale, we give sufficient conditions under which the generalized derivatives (of a given order) of the solution $u$ are continuous.

\begin{theorem}\label{1th5}
Let $0\leq p\in\mathbb{Z}$ and suppose that a function $u\in H^{m+1/2+}(\Omega)$ satisfies the condition of Theorem~$\ref{1th4}$ for a certain parameter $\varphi\in\mathrm{RO}$ such that $\sigma_{0}(\varphi)>m+1/2$ and
\begin{equation}\label{1f15}
\int_1^{\infty}t^{2p+n-1}\varphi^{-2}(t)dt<\infty.
\end{equation}
Then $u\in C^{p}(\Omega_{0}\cup\Gamma_{0})$.
\end{theorem}

\begin{remark}\label{1rem1}
Condition \eqref{1f15} is sharp in Theorem~\ref{1th5}. Namely, let $0\leq p\in\mathbb{Z}$ and $\varphi\in\mathrm{RO}$ with $\sigma_{0}(\varphi)>m+1/2$. Then it follows from the implication
\begin{equation}\label{implication}
\bigl(u\in H^{m+1/2+}(\Omega)\quad \mbox{and}\quad
(A,B)u\in\mathcal{H}^{\varphi\rho^{-2q}}_
{\mathrm{loc}}(\Omega_{0},\Gamma_{0})\bigr)\;\Rightarrow\;
u\in C^{p}(\Omega_{0}\cup\Gamma_{0})
\end{equation}
that $\varphi$ satisfies condition \eqref{1f15}.
\end{remark}

Theorem~\ref{1th5} implies a sufficient condition for the generalized solution $u$ to be classical.

\begin{corollary}\label{1prop6}
Let a function $u\in H^{m+1/2+}(\Omega)$ be a generalized solution to the elliptic boundary-value problem \eqref{1f1}, \eqref{1f2} in which
\begin{gather}\label{1f16}
f\in H^{\varphi_1\rho^{-2q}}_{\mathrm{loc}}(\Omega,\varnothing)\cap H^{\varphi_2\rho^{-2q}}(\Omega),\\
g_j\in H^{\varphi_2\rho^{-m_j-1/2}}(\Gamma),\quad\mbox{with}\quad
j=1,\ldots,q, \label{1f17}
\end{gather}
for some parameters $\varphi_1,\varphi_2\in\mathrm{RO}$ that satisfy the conditions $\sigma_0(\varphi_1)>m+1/2$, $\sigma_0(\varphi_2)>m+1/2$ and
\begin{gather}\label{f18}
\int_1^{\infty}t^{n+2q-1}\varphi_{1}^{-2}(t)dt<\infty,\\
\int_1^{\infty}t^{2m+n-1}\varphi_{2}^{-2}(t)dt<\infty. \label{f19}
\end{gather}
Then $u$ is a classical solution, i.e. $u\in C^{2q}(\Omega)\cap C^{m}(\overline{\Omega})$.
\end{corollary}

Note that if the solution $u$ is classical, then the left-hand sides of the problem \eqref{1f1}, \eqref{1f2} are calculated by means of classical derivatives and are continuous on $\Omega$ or $\Gamma$ respectively.

\section{Interpolation with function parameter}\label{1sec5}

An important property of the extended Sobolev scale consists in that this scale can be obtained by the interpolation with a function parameter between Sobolev inner product spaces. This interpolation will be a main tool in our proof of the key Theorem~\ref{1th1}. Therefore we recall the definition of this interpolation method and formulate its properties being used in the paper.

Note that the method of interpolation with a function parameter between normed spaces was introduced by Foia\c{s} and Lions in paper \cite{FoiasLions61}, where the case of Hilbert spaces was separately considered. As to Hilbert spaces, this method is a natural generalization of the classical interpolation method by Lions and S.~Krein (see, e.g., \cite[Chapter~1, Section 1 and 5]{LionsMagenes72}) to the case where a general enough function serves as an interpolation parameter instead of a number. Setting forth the interpolation with a function parameter between Hilbert spaces, we follow the monograph \cite[Section~1.1]{MikhailetsMurach14}. It is enough for our purposes to restrict ourselves to the case of separable Hilbert spaces.

Let $X:=[X_{0},X_{1}]$ be an ordered couple of separable complex Hilbert spaces
that the continuous and dense embedding $X_{1}\hookrightarrow X_{0}$ holds. We say
that this couple is admissible. As is known \cite[Chapter~1, Section~1]{LionsMagenes72}, for $X$ there exists an isometric isomorphism $J:X_{1}\leftrightarrow X_{0}$ such that $J$ is a self-adjoint positive operator on $X_{0}$ with the domain~$X_{1}$. The operator $J$ is called a generating operator for the couple~$X$. This operator is uniquely determined by $X$.

By $\mathcal{B}$ we denote the set of all Borel measurable functions
$\psi:(0,\infty)\rightarrow(0,\infty)$ such that $\psi$ is bounded on each compact
interval $[a,b]$, with $0<a<b<\infty$, and that $1/\psi$ is bounded on every set
$[r,\infty)$, with $r>0$.

Let $\psi\in\mathcal{B}$. Consider the operator $\psi(J)$, which is defined (and
positive) in $X_{0}$ as the Borel function $\psi$ of the positive self-adjoint operator $J$.
Denote by $[X_{0},X_{1}]_{\psi}$ or, simply, by $X_{\psi}$ the domain of the operator
$\psi(J)$ endowed with the inner product
$(w_1, w_2)_{X_\psi}:=(\psi(J)w_1,\psi(J)w_2)_{X_0}$ and the
corresponding norm $\|w\|_{X_\psi}=(w,w)_{X_\psi}^{1/2}$. The space $X_{\psi}$ is
Hilbert and separable. The continuous and dense embedding $X_\psi \hookrightarrow X_0$ holds true.

A function $\psi\in\mathcal{B}$ is called an interpolation parameter if the
following condition is fulfilled for all admissible couples $X=[X_{0},X_{1}]$ and
$Y=[Y_{0},Y_{1}]$ of Hilbert spaces and for an arbitrary linear mapping $T$ given on
$X_{0}$: if the restriction of $T$ to $X_{j}$ is a bounded operator
$T:X_{j}\rightarrow Y_{j}$ for each $j\in\{0,1\}$, then the restriction of $T$ to
$X_{\psi}$ is also a bounded operator $T:X_{\psi}\rightarrow Y_{\psi}$.
If $\psi$ is an interpolation parameter, then we say that the Hilbert space
$X_{\psi}$ is obtained by the interpolation of the couple $X=[X_{0},X_{1}]$
(or between $X_{0}$ and $X_{1}$) with the function parameter $\psi$. In this case,
the dense and continuous embeddings
$X_{1}\hookrightarrow X_{\psi}\hookrightarrow X_{0}$ hold.

Note that a function $\psi\in\mathcal{B}$ is an interpolation parameter if and only
if $\psi$ is pseudoconcave on a neighborhood of $+\infty$ (see \cite[Sect.~1.1.9]{MikhailetsMurach14}).
The latter condition means that there exists a concave
function $\psi_{1}:(b,\infty)\rightarrow(0,\infty)$, with $b\gg1$, such that both
functions $\psi/\psi_{1}$ and $\psi_{1}/\psi$ are bounded on $(b,\infty)$.
This fundamental result follows from Peetre's theorem \cite{Peetre66, Peetre68}
about description of all interpolation functions of positive order (see also
\cite[Sect.~5.4]{BerghLefstrem80}).

The above-mentioned interpolation property of the extended Sobolev is stated in the
following way.

\begin{proposition}\label{1prop1}
Assume that a function $\alpha\in\mathrm{RO}$ and numbers $s_{0},s_{1}\in\mathbb{R}$ satisfy the conditions $s_{0}<\sigma_{0}(\alpha)$ and $s_{1}>\sigma_{1}(\alpha)$. Put
\begin{equation}\label{1f20}
\psi(t):=
\begin{cases}
\;t^{{-s_0}/{(s_1-s_0)}}\,
\alpha\bigl(t^{1/{(s_1-s_0)}}\bigr)&\text{for}\quad t\geq1, \\
\;\alpha(1)&\text{for}\quad0<t<1.
\end{cases}
\end{equation}
Then $\psi\in\mathcal{B}$ is an interpolation parameter, and
\begin{equation*}
[H^{(s_0)}(G),H^{(s_1)}(G)]_{\psi}=H^{\alpha}(G)
\end{equation*}
with equality of norms if $G=\mathbb{R}^{n}$, and with equivalence of norms if $G=\Omega$ or $G=\Gamma$.
\end{proposition}

This proposition is proved by Mikhailets and Murach \cite[Theorem 5.1]{MikhailetsMurach15ResMath1} for $G=\Omega$ and in \cite[Theorems 2.19 and 2.22]{MikhailetsMurach14} for $G\in\{\mathbb{R}^{n},\Gamma\}$.

It is worthwhile to note that the extended Sobolev scale is closed with respect to the interpolation with a function parameter and coincides (up to equivalence of norms) with the class of all Hilbert
spaces that are interpolation spaces for the couples of the Sobolev
spaces $H^{(s_0)}(G)$ and $H^{(s_1)}(G)$ with $s_0,s_1\in\mathbb{R}$
and $s_0<s_1$. The latter property follows from Ovchinnikov's
theorem \cite[p.~511]{Ovchinnikov84}, which gives a description
of all Hilbert spaces that are interpolation spaces for an
arbitrarily chosen compatible couple of Hilbert spaces. Recall that the
property of a (Hilbert) space $H$ to be an interpolation space for
an admissible couple $X=[X_0,X_1]$ means the following: a)~the
continuous and dense embeddings $X_1\hookrightarrow H\hookrightarrow
X_0$ hold; b)~every linear operator which is bounded on each of the spaces
$X_0$ and $X_1$ should be bounded on $H$ as well.

Let us recall two general properties of the interpolation; they will be used in our proof of Theorem~\ref{1th1}.

\begin{proposition}\label{1prop2}
Let $X=[X_0,X_1]$ and  $Y=[Y_0,Y_1]$ be two admissible couples of  Hilbert spaces. Suppose that a linear mapping $T$ is given on $X_0$ and satisfies the following condition: the restrictions of $T$ to the spaces $X_0$ and $X_1$ are bounded and Fredholm operators $T:X_0\rightarrow Y_0$ and $T:X_1\rightarrow Y_1$ respectively, and these operators  have the same kernel and the same index. Then, for an arbitrary interpolation parameter $\psi\in\mathcal{B}$, the restriction of $T$ to $X_\psi$ is also a bounded and Fredholm operator $T:X_\psi\rightarrow Y_\psi$ with the same kernel and the same index; the range of this operator is equal to $Y_\psi\cap T(X_0)$.
\end{proposition}

\begin{proposition}\label{1prop3}
Let $[X_{0}^{(k)},X_{1}^{(k)}]$, with $k=1,\ldots,p$, be a finite
collection of admissible couples of Hilbert spaces. Then for every function $\psi\in\mathcal{B}$ we have
\begin{equation*}
\biggl[\,\bigoplus_{k=1}^{p}X_{0}^{(k)},\,
\bigoplus_{k=1}^{p}X_{1}^{(k)}\biggr]_{\psi}=\,
\bigoplus_{k=1}^{p}\bigl[X_{0}^{(k)},\,X_{1}^{(k)}\bigr]_{\psi}
\end{equation*}
with equality of norms.
\end{proposition}

These properties of the interpolation are expounded, e.g., in monograph \cite[Sections 1.1.5 and 1.1.7]{MikhailetsMurach14}.

As an application of this interpolation method, we obtain the following property of linear differential operators on the extended Sobolev scale.

\begin{proposition}\label{1prop4}
$\mathrm{(i)}$ Let $L$ be a linear differential operator on $\overline{\Omega}$ of order $\ell\geq0$ with coefficients from $C^{\infty}(\overline{\Omega})$. Then the mapping $u\mapsto Lu$, with $u\in C^{\infty}(\overline{\Omega})$, extends uniquely (by continuity) to the bounded linear operator
\begin{equation*}
L:\,H^\alpha(\Omega)\rightarrow H^{\alpha\rho^{-\ell}}(\Omega)
\end{equation*}
for every $\alpha\in\mathrm{RO}$.

$\mathrm{(ii)}$ Let $K$ be a boundary  linear differential operator on $\Gamma$ of order $k\geq\nobreak0$
with coefficients from $C^{\infty}(\Gamma)$. Then the mapping $u\mapsto Ku$, with $u\in C^{\infty}(\overline{\Omega})$, extends uniquely (by continuity) to the bounded linear operator
\begin{equation*}
K:\,H^\alpha(\Omega)\rightarrow H^{\alpha\rho^{-k-1/2}}(\Gamma)
\end{equation*}
for every function parameter $\alpha\in\mathrm{RO}$ such that $\sigma_{0}(\alpha)>k+1/2$.
\end{proposition}

Proposition \ref{1prop4} is well-known in the case of Sobolev spaces. For the extended Sobolev scale, this
proposition follows directly from the Sobolev case due to Proposition~\ref{1prop1}.

\section{Proofs}\label{1sec6}

In this section, we will prove Theorems \ref{1th1}--\ref{1th5}, Remark~\ref{1rem1}, and Corollary~\ref{1prop6}.

\medskip

\noindent\emph{Proof of Theorem~$\ref{1th1}$.} It follows directly from Proposition~\ref{1prop4} that the mapping $C^{\infty}(\overline{\Omega})\ni u\mapsto(Au,Bu)$ extends uniquely to the bounded linear operator \eqref{1f8}. We first discuss Theorem~\ref{1th1} in the Sobolev case and then deduce it in the general situation with the help of the interpolation with a function parameter between Sobolev spaces.

Consider Theorem~\ref{1th1} in the Sobolev case where $\varphi(t)\equiv t^{s}$ and $s>m+1/2$. The central and well-known
fact of the theory of elliptic boundary-value problems in Sobolev spaces consists in the following: the operator
\eqref{1f8} is Fredholm for every real $s\geq2q$, its kernel coincides with $N$, and its index does not depend on $s$
(see, e.g., survey \cite[Section~2.4]{Agranovich97}). It is less known that this fact holds true for arbitrary real
$s>m+1/2$ (see, e.g., \cite[Chapter~III, Section~2.2]{Egorov84}). The proof of the remaining part of
Theorem~\ref{1th1} -- concerning the range of \eqref{1f8} -- is given in monographs \cite[Theorem 3.4.1 and Section 4.4.1]{KozlovMazyaRossmann97} and \cite[Theorem 2.4.1]{Roitberg99} for every real $s\geq2q$. Therefore it remains to prove this part in the case where $m+1/2<s<2q$.

Let $\Lambda_{s}$ denote the Fredholm bounded operator \eqref{1f8} in the Sobolev case where $\varphi(t)\equiv t^{s}$ and $s>m+1/2$. As usual, $\Lambda_{s}^{*}$ stands for the adjoint operator to \eqref{1f8}. If $m+1/2<s_{0}<s_{1}$, then $\Lambda_{s_{0}}^{\ast}$ is a restriction of $\Lambda_{s_{1}}^{\ast}$ and therefore $\ker\Lambda_{s_{0}}^{\ast}\subseteq\ker\Lambda_{s_{1}}^{\ast}$. But  $\ker\Lambda_{s_{0}}^{\ast}$ and $\ker\Lambda_{s_{1}}^{\ast}$ have the same finite dimension. Indeed, the latter is equal to the co-dimension of the range of $\Lambda_{s}$, and this codimension is finite and does not depend on $s$, as we have mentioned. Thus, $\ker\Lambda_{s_{0}}^{\ast}=\ker\Lambda_{s_{1}}^{\ast}$. According to  \cite[Theorem 3.4.2]{KozlovMazyaRossmann97}, we have the equality $\ker\Lambda_{s}^{\ast}=N_{\ast}$ if $s\geq 2q$. Here, every vector $(v,h_{1},\ldots,h_{q})\in N_{\ast}$ is interpreted as the continuous linear functional
\begin{equation*}
\Upsilon:(f,g_{1},\ldots,g_{q})\mapsto
(f,v)_\Omega+\sum_{j=1}^{q}(g_{j},h_{j})_{\Gamma}
\end{equation*}
on the space
\begin{equation*}
\mathcal{H}^{(s-2q)}(\Omega,\Gamma):=H^{(s-2q)}(\Omega)
\oplus\bigoplus_{j=1}^{q}H^{(s-m_j-1/2)}(\Gamma).
\end{equation*}
Hence, $\ker\Lambda_{s}^{\ast}=N_{\ast}$ whenever $s>m+1/2$. Therefore, $\Upsilon$ extends uniquely (by continuity) to a linear continuous functional on $\mathcal{H}^{(s-2q)}(\Omega,\Gamma)$ in the case where $m+1/2<s<2q$. Specifically, the form $(f,v)_\Omega$ is well defined for every $f\in H^{(s-2q)}(\Omega)$ via the passage to the limit in this case. Now, since the operator $\Lambda_{s}$ is Fredholm whenever $s>m+1/2$, then its range consists of all vectors $(f,g_{1},\ldots,g_{q})\in\mathcal{H}^{(s-2q)}(\Omega,\Gamma)$ such that $\Upsilon(f,g_{1},\ldots,g_{q})=\nobreak0$ for all $\Upsilon\in\ker\Lambda_{s}^{\ast}$. This condition becomes \eqref{1f9} in view of the equality $\ker\Lambda_{s}^{\ast}=N_{\ast}$.
Besides, the index of $\Lambda_{s}$ is $\dim N-\dim N^{\ast}$ because $\dim\ker\Lambda_{s}^{\ast}$ is equal to the co-dimension of the range of $\Lambda_{s}$.

Thus, Theorem~\ref{1th1} is justified in the Sobolev case. We now  deduce this theorem in the general situation from this case with the help of the interpolation.

We arbitrarily choose $\varphi\in\mathrm{RO}$ such that $\sigma_0(\varphi)>m+1/2$. Let real numbers $l_{0}$ and $l_{1}$ satisfy the conditions $m+1/2<l_{0}<\sigma_{0}(\varphi)$ and $\sigma_{1}(\varphi)< \nobreak l_{1}$. We have the Fredholm bounded operators
\begin{equation}\label{1f21}
\Lambda_{l_{i}}=(A,B):\,H^{(l_{i})}(\Omega)\rightarrow
\mathcal{H}^{(l_{i}-2q)}(\Omega,\Gamma),\quad
i\in\{0,1\}.
\end{equation}
They have the common kernel $N$ and the same index $\dim N-\dim N_{*}$. Besides,
\begin{gather}\label{1f22}
(A,B)\bigl(H^{(l_{i})}(\Omega)\bigr)=
\bigl\{(f,g)\in\mathcal{H}^{(l_{i}-2q)}(\Omega,\Gamma):\,
\mbox{\eqref{1f9} is true}\bigr\}.
\end{gather}

We define an interpolation parameter $\psi$ by formula \eqref{1f20} with $\alpha:=\varphi$, $s_{0}:=l_{0}$, and $s_{1}:=\nobreak l_{1}$. Applying the interpolation with the function parameter $\psi$ to \eqref{1f21}, we conclude by Proposition~\ref{1prop2} that the restriction of $\Lambda_{l_{0}}$ to the interpolation space $[H^{(l_{0})}(\Omega),H^{(l_{1})}(\Omega)]_{\psi}$ is a bounded and Fredholm  operator
\begin{equation}\label{1f23}
(A,B):\,\bigl[H^{(l_{0})}(\Omega),H^{(l_{1})}(\Omega)\bigr]_{\psi}\to
\bigl[\mathcal{H}^{(l_{0}-2q)}(\Omega,\Gamma),
\mathcal{H}^{(l_{1}-2q)}(\Omega,\Gamma)\bigr]_{\psi}.
\end{equation}
Since $C^{\infty}(\overline{\Omega})$ is dense in this space, operator \eqref{1f23} is an extension by continuity of the mapping $C^{\infty}(\overline{\Omega})\ni u\mapsto(Au,Bu)$. This operator coincides with \eqref{1f8} because
\begin{gather*}
\bigl[H^{(l_{0})}(\Omega),H^{(l_{1})}(\Omega)\bigr]_{\psi}=
H^{\varphi}(\Omega),\\
\bigl[\mathcal{H}^{(l_{0}-2q)}(\Omega,\Gamma),
\mathcal{H}^{(l_{1}-2q)}(\Omega,\Gamma)\bigr]_{\psi}=
\mathcal{H}^{\varphi\rho^{-2q}}(\Omega,\Gamma) \label{1f24b}
\end{gather*}
with equivalence of norms.

Indeed, the first formula is due to Proposition~\ref{1prop1}, in which $\alpha:=\varphi$, $s_{0}:=l_{0}$, $s_{1}:=l_{1}$, and $G:=\Omega$. The second formula follows from Propositions \ref{1prop3} and \ref{1prop1}; namely,
\begin{align*}
&\bigl[\mathcal{H}^{(l_{0}-2q)}(\Omega,\Gamma),
\mathcal{H}^{(l_{1}-2q)}(\Omega,\Gamma)\bigr]_{\psi}\\
& \quad  =\bigl[H^{(l_{0}-2q)}(\Omega),H^{(l_{1}-2q)}(\Omega)\bigr]_{\psi}\oplus
\bigoplus_{j=1}^{q}\bigl[H^{(l_{0}-m_{j}-1/2)}(\Gamma),
H^{(l_{1}-m_{j}-1/2)}(\Gamma)\bigr]_{\psi}\\
& \quad =H^{\varphi\rho^{-2q}}(\Omega)\oplus\bigoplus_{j=1}^{q}
H^{\varphi\rho^{-m_j-1/2}}(\Gamma)=
\mathcal{H}^{\varphi\rho^{-2q}}(\Omega,\Gamma).
\end{align*}
Note that in Proposition~\ref{1prop1} we have first taken
$$
\alpha:=\varphi\rho^{-2q},\quad s_{0}:=l_{0}-2q,\quad s_{1}:=l_{1}-2q,
\quad G:=\Omega
$$
and then
\begin{equation*}
\alpha:=\varphi\rho^{-m_{j}-1/2},\quad
s_{0}:=l_{0}-m_{j}-1/2,\quad s_{1}:=l_{1}-m_{j}-1/2,\quad G:=\Gamma.
\end{equation*}
Formula \eqref{1f20} gives the same function $\psi$ in both of these cases.

According to Proposition~\ref{1prop2}, the kernel and index of the operator \eqref{1f8} coincide with the common kernel  $N$ and index $\dim N-\dim N_{*}$ of the operators \eqref{1f21}. Besides, the range of \eqref{1f8} equals
\begin{equation}\label{domain}
\mathcal{H}^{\varphi\rho^{-2q}}(\Omega,\Gamma)\cap
(A,B)\bigl(H^{(l_{0})}(\Omega)\bigr)=
\bigl\{(f,g)\in\mathcal{H}^{\varphi\rho^{-2q}}(\Omega,\Gamma):\,
\mbox{\eqref{1f9} is true}\bigr\}.
\end{equation}
Here, we use the equality \eqref{1f22} with $i:=0$ and the embedding $$
\mathcal{H}^{\varphi\rho^{-2q}}(\Omega,\Gamma)\hookrightarrow
\mathcal{H}^{(l_{0}-2q)}(\Omega,\Gamma).
$$
Thus, the operator \eqref{1f8} has all the properties stated in  Theorem~\ref{1th1}.

This theorem is proved.

\medskip

\noindent\emph{Proof of Theorem $\ref{1th2}$.} According to Theorem~\ref{1th1}, the bounded linear operator \eqref{1f12} is a bijection. Hence it is an isomorphism due to the Banach theorem on inverse operator. Theorem~\ref{1th2} is proved.

\medskip

\noindent\emph{Proof of Theorem $\ref{1th3}$.} At first we will prove  the estimate \eqref{1f14} in the case where $\chi=\eta=1$. According to  \eqref{1f10}, we write an arbitrary function $u\in H^{\varphi}(\Omega)$ in the form $u=u_{0}+u_{1}$ with $u_{0}:=(1-P)u\in N$ and $u_{1}:=Pu\in P(H^{\varphi}(\Omega))$. By virtue of Theorem~\ref{1th2}, we have
\begin{equation}\label{t}
\|u_1\|_{H^{\varphi}(\Omega)}\leq
c_{1}\|(A,B)u_1\|_{\mathcal{H}^{\varphi\rho^{-2q}}(\Omega,\Gamma)}=
c_{1}\|(A,B)u\|_{\mathcal{H}^{\varphi\rho^{-2q}}(\Omega,\Gamma)}.
\end{equation}
Here, $c_1$ is the norm of the inverse operator to the isomorphism \eqref{1f12} (note also that $(A,B)u_{1}=(A,B)u$). Since the space $N$ is finite-dimensional, all norms are equivalent on $N$, specifically, the norms in $H^{\varphi}(\Omega)$ and $H^{\varphi\rho^{-1}}(\Omega)$. Hence, in view of $u_{0}\in N$ and \eqref{t}, we obtain the following inequalities:
\begin{align*}
\|u_0\|_{H^{\varphi}(\Omega)}&\leq c_{0}\|u_0\|_{H^{\varphi\rho^{-1}}(\Omega)}\leq
c_{0}\|u\|_{H^{\varphi\rho^{-1}}(\Omega)}+
c_{0}\|u_1\|_{H^{\varphi\rho^{-1}}(\Omega)}\\
&\leq c_{0}\|u\|_{H^{\varphi\rho^{-1}}(\Omega)}+
c_{0}\|u_1\|_{H^{\varphi}(\Omega)}\\
&\leq c_{0}\|u\|_{H^{\varphi\rho^{-1}}(\Omega)}+
c_{0}c_{1}\|(A,B)u\|_{\mathcal{H}^{\varphi\rho^{-2q}}(\Omega,\Gamma)}.
\end{align*}
Here, $c_0$ is a certain positive number that does not depend on $u$. So,
\begin{equation}\label{tt}
\|u_0\|_{H^{\varphi}(\Omega)}\leq c_{0}c_{1}\|(A,B)u\|_
{\mathcal{H}^{\varphi\rho^{-2q}}(\Omega,\Gamma)}+
c_{0}\|u\|_{H^{\varphi\rho^{-1}}(\Omega)}.
\end{equation}
The inequalities \eqref{t} and \eqref{tt} gives the estimate \eqref{1f14} in the case of $\chi=\eta=1$; namely,
\begin{equation}\label{global-estimate}
\|u\|_{H^{\varphi}(\Omega)}\leq c_3
\bigl(\|(A, B)u\|_{\mathcal{H}^{\varphi\rho^{-2q}}(\Omega,\Gamma)}+
\|u\|_{H^{\varphi\rho^{-1}}(\Omega)}\bigr)
\end{equation}
for arbitrary $u\in H^{\varphi}(\Omega)$. Here, the positive number $c_3:=c_0+c_1+c_0c_1$ does not depend on $u$.

We will now deduce Theorem~\ref{1th3} in the general situation from  \eqref{global-estimate}. Choose a function $u\in H^{\varphi}(\Omega)$ arbitrarily. In our reasoning, $c_3$, $c_4$, and $c_{5}$ denote some positive numbers that are independent of~$u$.

Consider the estimate \eqref{global-estimate} with $\chi u\in H^{\varphi}(\Omega)$ instead of $u$, namely,
\begin{equation}\label{1f29}
\|\chi u\|_{H^{\varphi}(\Omega)}\leq c_3\,\bigl(\|(A, B)(\chi u)\|_{\mathcal{H}^{\varphi\rho^{-2q}}(\Omega,\Gamma)}+
\|\chi u\|_{H^{\varphi\rho^{-1}}(\Omega)}\bigr).
\end{equation}
Rearranging the operator of multiplication by $\chi$ with the differential operators $A$ and $B_1,\ldots,B_q$, we write
\begin{align*}
(A,B)(\chi u)&=(A,B)(\chi\eta u)=\chi(A,B)(\eta u)+(A',B')(\eta u)\\
&=\chi(A,B)u+(A',B')(\eta u).
\end{align*}
Here, $A'$ is a linear differential operator on $\overline{\Omega}$ of order $\mathrm{ord}\,A'\leq 2q-1$, and $B':=(B_{1}',\ldots,B_{q}')$ is a collection of boundary linear differential operators on $\Gamma$ of order $\mathrm{ord}\,B_{j}'\leq m_{j}-1$ for each $j\in\{1,\ldots,q\}$. All the coefficients of $A'$ and $B_{j}'$ belong respectively to $C^{\infty}(\overline{\Omega})$ and $C^{\infty}(\Gamma)$. Thus,
\begin{equation}\label{1f30}
(A,B)(\chi u)=\chi(A,B)u+(A',B')(\eta u),
\end{equation}
where, by Proposition \ref{1prop4},
\begin{equation}\label{1f30b}
\|(A',B')(\eta u)\|_{\mathcal{H}^{\varphi\rho^{-2q}}(\Omega,\Gamma)}
\leq c_{4}\|\eta u\|_{H^{\varphi\rho^{-1}}(\Omega)}.
\end{equation}
Now, according to \eqref{1f29}--\eqref{1f30b}, we obtain
\begin{equation}\label{1f30c}
\begin{aligned}
\|\chi u\|_{H^{\varphi}(\Omega)}&\leq c_3\,\bigl(
\|\chi(A, B)u\|_{\mathcal{H}^{\varphi\rho^{-2q}}(\Omega,\Gamma)}\\
&\quad+\|(A',B')(\eta u)\|_{\mathcal{H}^{\varphi\rho^{-2q}}(\Omega,\Gamma)}+
\|\chi u\|_{H^{\varphi\rho^{-1}}(\Omega)}\bigr)\\
&\leq c_3\|\chi(A, B)u\|_
{\mathcal{H}^{\varphi\rho^{-2q}}(\Omega,\Gamma)}+
c_{3}c_{4}\|\eta u\|_{H^{\varphi\rho^{-1}}(\Omega)}+
c_3\|\chi u\|_{H^{\varphi\rho^{-1}}(\Omega)}.
\end{aligned}
\end{equation}
Note that, by Proposition \ref{1prop4},
\begin{equation}\label{1f30d}
\|\chi u\|_{H^{\varphi\rho^{-1}}(\Omega)}=
\|\chi\eta u\|_{H^{\varphi\rho^{-1}}(\Omega)}\leq
c_{5}\|\eta u\|_{H^{\varphi\rho^{-1}}(\Omega)}.
\end{equation}
Formulas \eqref{1f30c} and \eqref{1f30d} immediately give the required estimate \eqref{1f14}, in which $c:=c_{3}(1+c_{4}+c_{5})$.

Theorem \ref{1th3} is proved.

\medskip

\noindent\emph{Proof of Theorem $\ref{1th4}$.} We will first prove this theorem in the case where $\Omega_{0}=\Omega$ and $\Gamma_{0}=\Gamma$. By the condition, $u\in H^{(s)}(\Omega)$ for some real $s\in(m+1/2,\sigma_{0}(\varphi))$, and
$$
(f,g):=(f,g_{1},\ldots,g_{q})=
(A,B)u\in\mathcal{H}^{\varphi\rho^{-2q}}(\Omega,\Gamma).
$$
Therefore, owing to Theorem \ref{1th1} applied also in the Sobolev case, the vector $(f,g)$ satisfies condition \eqref{1f9}, and hence $(A,B)v=(f,g)$ for some function $v\in H^{\varphi}(\Omega)$. Thus, $(A,B)(u-v)=0$, with $u-v\in H^{(s)}(\Omega)$. Then $w:=u-v\in N\subset C^{\infty}(\overline{\Omega})$ due to Theorem~\ref{1th1}. Hence, $u=v+w\in H^{\varphi}(\Omega)$. Theorem~\ref{1th4} is proved in the case considered.

Let us now prove this theorem in the general situation. We arbitrarily choose an open set $V_{1}\subset\mathbb{R}^{n}$ such that $\overline{V_{1}}\subset V$ and $\Omega\cap V_{1}\neq\varnothing$. Put  $\Omega_1:=\Omega\cap V_1$ and $\Gamma_{1}:=\Gamma\cap V_1$. Let us prove that $u\in H^{\varphi}_{\mathrm{loc}}(\Omega_{1},\Gamma_{1})$.

We choose functions $\chi,\eta\in C^{\infty}(\overline{\Omega})$ so that their supports lie in $\Omega_{0}\cup\Gamma_{0}$, $\chi=1$ in a neighborhood of $\Omega_{1}\cup\Gamma_{1}$, and $\eta=1$ on   $\mathrm{supp}\,\chi$. By the condition, $u\in H^{(s)}(\Omega)$ for some real $s\in(m+1/2,\sigma_{0}(\varphi))$ and $(A,B)u=(f,g)$ with $(f,g)\in
\mathcal{H}^{\varphi\rho^{-2q}}_{\mathrm{loc}}(\Omega_{0},\Gamma_{0})$. Therefore
\begin{equation*}
(A,B)(\chi u)=\eta(A,B)(\chi u)=\eta(f,g)-\eta(A,B)((1-\chi)u).
\end{equation*}
Using the projector $P_{\ast}$, we write
\begin{equation*}
(A,B)(\chi u)=P_{\ast}(\eta(f,g))+F,
\end{equation*}
where
\begin{gather*}
P_{\ast}(\eta(f,g))\in P_{\ast}\bigl(\mathcal{H}^{\varphi\rho^{-2q}}(\Omega,\Gamma)\bigr),\\
F:=(1-P_{\ast})(\eta(f,g))-\eta(A,B)((1-\chi)u)\in
P_{\ast}\bigl(\mathcal{H}^{\varrho^{s-2q}}(\Omega,\Gamma)\bigr).
\end{gather*}
Hence, by Theorem \ref{1th2}, there exist functions $u_1\in H^{\varphi}(\Omega)$ and $u_2\in H^{(s)}(\Omega)$ such that $(A,B)u_1=P_{\ast}(\eta(f,g))$ and $(A,B)u_2=F$. Then $(A,B)(\chi u-u_1-u_2)=0$, which implies
\begin{equation*}
w:=\chi u-u_1-u_2\in N\subset C^{\infty}(\overline{\Omega})
\end{equation*}
in view of Theorem~\ref{1th1}. We note that $F\in\mathcal{H}^{\varrho^{l-2q}}_{\mathrm{loc}}(\Omega_{1},\Gamma_{1})$ for every $l>\sigma_{1}(\varphi)$ because $(1-P_{\ast})(\eta(f,g))\in N_{\ast}$ and $\eta(A,B)((1-\chi)u)=0$ on $\Omega_{1}\cup\Gamma_{1}$. Therefore
\begin{equation*}
u_2\in H_{\mathrm{loc}}^{\varrho^l}(\Omega_{1},\Gamma_{1})\subset
H^{\varphi}_{\mathrm{loc}}(\Omega_{1},\Gamma_{1})
\end{equation*}
according to the theory of elliptic boundary-value problems in Sobolev spaces (see, e.g., \cite[Theorem~5.2]{MikhailetsMurach07UMJ5}). We conclude now that
\begin{equation*}
\chi u=u_1+u_2+w\in H^{\varphi}_{\mathrm{loc}}(\Omega_{1},\Gamma_{1}).
\end{equation*}
Then $\zeta u=\zeta\chi u\in H^{\varphi}(\Omega)$ for an arbitrary function $\zeta\in C^{\infty}(\overline{\Omega})$ with $\mathrm{supp}\,\zeta\subset\Omega_{1}\cup\Gamma_{1}$; i.e., $u\in H^{\varphi}_{\mathrm{loc}}(\Omega_{1},\Gamma_{1})$. The latter inclusion implies that
$u\in H^{\varphi}_{\mathrm{loc}}(\Omega_{0},\Gamma_{0})$ due to our choice of $V_{1}$.

Theorem \ref{1th4} is proved.

\medskip

\noindent\emph{Proof of Theorem $\ref{1th5}$.} According to Theorem~\ref{1th4}, the inclusion $u\in H^{\varphi}_{\mathrm{loc}}(\Omega_{0},\Gamma_{0})$ holds.  We arbitrarily choose a point $x\in\Omega_{0}\cup\Gamma_{0}$. Let a function $\chi\in C^{\infty}(\overline{\Omega})$ satisfy the conditions $\mathrm{supp}\,\chi\subset\Omega_0\cup\Gamma_0$ and $\chi=1$ in some
neighborhood $V(x)$ of~$x$. Owing to \eqref{v}, we get
$\chi u\in H^{\varphi}(\Omega)\subset C^{p}(\overline{\Omega})$.
Then $u\in C^{p}(V(x))$. Hence, $u\in
C^{p}(\Omega_{0}\cup\Gamma_{0})$ due to the arbitrariness of $x$.
Theorem~\ref{1th5} is proved.

\medskip

\noindent\emph{Proof of Remark $\ref{1rem1}$.} Suppose that the implication \eqref{implication} holds true. Let $V$ be an open ball such that $\overline{V}\subset\Omega_{0}$, and choose a function $v\in H^{\varphi}(V)$ arbitrarily. Note that $v=u\!\upharpoonright\!V$ for certain $u\in H^{\varphi}(\Omega)$. Besides, $(A,B)u\in\mathcal{H}^{\varphi\rho^{-2q}}(\Omega,\Gamma)$. Now, according to \eqref{implication}, we conclude that $u\in C^{p}(\Omega_{0}\cup\Gamma_{0})$. Hence, $v\in C^{p}(\overline{V})$. Thus, $H^{\varphi}(V)\subset C^{p}(\overline{V})$, which implies
condition \eqref{1f15} in view of \eqref{v}. Remark~\ref{1rem1} is proved.

\medskip

\noindent\emph{Proof of Corollary $\ref{1prop6}$.} By virtue of Theorem~\ref{1th5} with $p:=2q$, $\varphi:=\varphi_1$, $\Omega_{0}:=\Omega$, and $\Gamma_{0}:=\varnothing$, it follows from conditions \eqref{1f16} and  \eqref{f18} that $u\in C^{2q}(\Omega)$. Besides, by the same theorem with $p:=m$, $\varphi:=\varphi_2$, $\Omega_{0}:=\Omega$, and $\Gamma_{0}:=\Gamma$, it follows from conditions \eqref{1f16}, \eqref{1f17}, and \eqref{f19} that $u\in C^{m}(\overline{\Omega})$. Thus, $u$ is a classical solution to the boundary-value problem  \eqref{1f1}, \eqref{1f2}. Corollary~\ref{1prop6} is proved.


\end{document}